# Destruction of the Resident Enterprise in the Special Economic Zone with Sanctions

Sergei Masaev
Institute of Oil and Gas
Siberian Federal University
Krasnoyarsk, Russia
Smasaev@sfu-kras.ru

*Abstract*—The activity of a special economic zone is defined by a dynamic equation, taking into account the individual strategies of residents. At a given point in time, in respect to the resident enterprise of a special economic zone, a regime is introduced that limits the flow of resources by 80% (sanctions), forming an integral indicator for a comprehensive assessment of the impact of sanctions on the enterprise. On the basis of the dynamic equation, an estimate of the economic damage for the potential SEZ of the Krasnoyarsk Territory is given.

## I. Introduction

Economic sanctions have existed for hundreds of years. In fact, the first economic sanctions, recorded by history, dates back to 423 BC. Athens, as the dominant economy of Hellas, did not allow traders from the Megara region to its ports and markets. For more than 20 centuries, the use of economic sanctions has been the prerogative of the dominant economy in the territory in question. Sanctions were introduced and are introduced by all leading countries, such as: the United States, the European Union countries, the USSR, Russia, Japan, China, Ukraine, Belarus and Israel. The United Kingdom achieved particular unsurpassed success in the conduct of trade restrictions in 1888. With 2% of the world's population, it controlled 54% of the turnover of the world economy.

Russia, in its documented history, faced the definition of "economic sanctions" in 1548, while implementing the plan to recruit the best German artisans to work in Russia by order of Ivan the Terrible. From then until today, there has not been a single period, with the growing influence of Russia, in which no sanctions were introduced against it that preceded military conflicts. The introduction of new "economic sanctions" by modern Russia has been felt since 2008. In the 21st century, the main economic sanctions by the United States were being introduced as an economic payment for Russia's political action in the Ossetian-Georgian conflict of August 2008 and the annexation of the Crimea in 2014. It is important to note that the United States since 1931 has the most developed global network of special economic zones allowing it to dictate conditions to China and significantly influence the economy of the European Union.

## II. Formulation of the problem

Obvious and implicit economic sanctions will exist in any economic relations between entities, regions, countries and political unions, therefore it is necessary to create a tool to manage this activity. The difficulties lie in the absence of indicators reflecting: a dynamic change in the form of economic interaction between large business entities, the influence of the external environment and sanctions. Traditionally in this area a lot of economic work of is done in a general nature [1-3]. In control theory, the beginning of the development of economic and mathematical models is associated with the names of V.V. Leontiev, V.S. Nemchinova, L.V. Kantorovich [4-6]. Further successes in the development of modeling methods were achieved in the 50s – 60s of the last century in the works of R. Bellman, K.A. Bagrinovskogo, A.G. Granberg, A.G. Aganbegyan, V.L. Makarova, Yu.P. Ivanilov and other authors. Since the late 1970s, a great contribution to the development of dynamic mathematical modeling of production processes, taking into account the raw material and production base, has been made by V.F. Krotov, E. Deming, T.K. Sirazetdinov, P.K. Semenov, Yu.P. Shurgin, C.K. Dzhaksybayeva, I.S. Ivanenko, V.V. Rodionov, A.A. Afanasyev [7-9].

Separately, it is worth noting Burkov V.N., who created, around the turn at the junction of the 21st century, a theory of active systems based on mathematical modeling (system analysis, game theory, decision theory, operations research), management decision-making procedures in organizational systems (economic agents ) [10-11] and developments in the following areas: mechanisms for controlling multi-agent active systems (D.A. Novikov, N.A. Korgin, A.V. Schepkin); theory of evolutionary development (V.V. Tsyganov); multi-level active systems and problems of the synthesis of structures, including network structures (M.V. Gubko, S.P. Mishin); information management mechanisms that take into account the reflexive behavior of agents (D.A. Novikov, A.G. Chkhartishvili); project management mechanisms (VN Burkov, A.V. Tsvetkov); mechanisms of innovation development (R.M. Nizhegorodtsev) [12].

In addition to the classical methods of the theory of management of the original approach to enterprise management is the method of integral indicators. In 2009, this method was used to analyze the state of the «Novy Gorod» enterprise in the context of the 2008 global financial crisis. Observing the non-stationary (dynamically changing set of controlled parameters and resources) state of the mathematical model of the observed enterprise allowed for the prediction and repairment of the functions performed by the enterprise caused by the global crisis of 2009 [13-17].

The purpose of this work is *to use the method of integrated indicators to predict the state of the enterprise resident of the special economic zone (SEZ) of the Krasnoyarsk Territory in the sanctions regime.*

### III. RESEARCH METHODOLOGY

A special economic zone (SEZ) is represented by a dynamic equation, as a dynamic system with external constraints from 8 resident enterprises:

$$y(t) = A(t)x(t) + B(t)u(t) + v(t). \quad (1)$$

$C = \{c_1, \ldots, c_i\}$ - list of companies, $I$ - their total number.
$T = \{t: t = 1, \ldots, T_{max}\}$ - set of time, many time points (months).
$x(t) = [x_1^1(t), x_2^2(t), \ldots, x_N^i(t)]^T$ - $N$ - economic parameter vector $i$ enterprise, where $x_N^i(t)$ value $n$ cost/income $i$ enterprise at the time $t$ of subspace $X^i$ of space $\Omega$.

There is some $g$ regulation for planning the allocation of available resources to functions $x_i$ for an enterprise based on the values of past periods $x_i^*$, then $x = g(x_i^*)$ with criterion

$$X_i(t) = \sum_{t=1}^{T} \sum_{i=1}^{n} TX^i \to max. \quad (2)$$

The planning function $x^*(t)$ is performed according to past economic indicators with a lag of $l$ periods $x(t) = p(x(t-l), \varepsilon(t))$, if $\varepsilon(t)$ is an error equal to 0, then the plan is equal to the fact $x^*(t) = x(t)$. A detailed description of the management of the enterprise, according to six strategy options, is presented in the work by reference [15].

$u(t) = [u_1(t), u_2(t), \ldots, u_M(t)]^T$ - $M$ is control vector, where $u_i(t)$ is the above control actions from the state in the form of changes in tax interest rates, industry subsidies and benefits aimed at increasing GRP (gross regional product) at time $t$:

$u_1$ - the rate of income tax in the budget of the subject of the Russian Federation,
$u_2$ - the rate of income tax in the consolidated budget of the Russian Federation,
$u_3$ - transport tax,
$u_4$ - property tax,
$u_5$ - the cost of electricity,
$u_6$ - tariffs for transportation of products,
$u_7$ - the average value of contributions to the pension fund of the Russian Federation (PF), the territorial fund for compulsory medical insurance (MHIF) of the Krasnoyarsk Territory, the federal fund for compulsory medical insurance (FFOMS), the social insurance fund of the Russian Federation (FSS),
$u_8$ - is the rental price of the land plot (land rent),
$u_9$ - the rental price of the forest plot (forest rent).

$y(t) = [y_1(t), y_2(t), \ldots, y_K(t)]^T$ - $K$ is the vector of observations, where $y_i(t)$ is the observed values of GRP (gross regional product) at time $t$.
$v(t) = [v_1(t), v_2(t), \ldots, v_K(t)]^T$ is a disturbance acting on $x(t)$ or other known factors that we cannot influence $v_{10}$ is the dollar rate (70 rubles / dollar is assumed in the model).

Other parameters:
$v_1$ - prices for resources,
$v_2$ - outstripping salary growth (in the model + 4% annually),
$v_3$ - investments from the owner,
$v_4$ - technological innovations,
$v_5$ - the movement of material flows,
$v_6$ - measures to improve the logistics of the project,
$v_7$ - labor resources,
$v_8$ - technology prices,
$v_9$ - inflation (in the model 4% annually).

$A_{N \times N}(t)$ - the matrix, which determines the speed of enterprise development, due to the use of financial resources and the influence of other enterprises. $a_{ij}$ - the degree of influence of economic parameters on each other $x_i$ on $x_j$.

$B_{M \times N}(t)$ - the matrix that determines the development of the enterprise, with changes in tax rates, subsidies and benefits. $b_{ij}$ - the degree of influence of the controlling influence $u_j(t)$ on the development of the $j$-th economic parameter of the enterprise $x_j(t)$.

$H_{K \times N}(t)$ - matrix of observations of GRP, which allows one to obtain an estimate of GRP $y_j(t)$ at the actual level $x_j(t)$.

To calculate the integral indicator $G_i$, it is necessary to represent the enterprise as a dynamic system at the moment with the considered parameters $x_n^i(t+l)$. Then we have a matrix $X^i$

$$X^i = \begin{pmatrix} x_1^i(t) & x_1^i(t+1) & x_1^i(t+l) & \cdots & x_1^i(T_{max}) \\ x_2^i(t) & x_2^i(t+1) & x_2^i(t+l) & \cdots & x_2^i(T_{max}) \\ \vdots & \vdots & \vdots & \ddots & \vdots \\ x_n^i(t) & x_n^i(t+1) & x_n^i(t+l) & \cdots & x_n^i(T_{max}) \end{pmatrix}. \quad (2)$$

Next, we calculate the coefficients of mutual correlation between the values of the parameters characterizing the state of the system for the entire planning period. To do this, it is necessary to calculate the correlation matrix $R^i(T_{max})$ at time points for all $t$, where $l$ is the time lag.

$$r_{t,t+l}(t) = \frac{\sum_{j=1}^{n}\left(x_j^i(t) - \bar{x}(t)\right) \cdot \left(x_j^i(t+l) - \bar{x}(t+l)\right)}{\sqrt{\sum_{j=1}^{n}\left(x_j^i(t) - \bar{x}(t)\right)^2} \cdot \sqrt{\sum_{j=1}^{n}\left(x_j^i(t+l) - \bar{x}(t+l)\right)^2}}, \quad (3)$$

$t$ - time points; $r_{t,t+l}(t)$ - correlation coefficients of variables $\sum_{j=1}^{n} x_j^i(t)$ and $\sum_{j=1}^{n} x_j^i(t+l)$ for all $x_n^i(t)$.

Form the matrix:

$$R^i(T_{max}) = \begin{pmatrix} r^i_{1,1}(1) & r^i_{1,2}(1) & r^i_{1,t+l}(1) & \cdots & r^i_{1,T_{max}}(1) \\ r^i_{t,1}(2) & r^i_{t,2}(2) & r^i_{t,t+l}(2) & \cdots & r^i_{t,T_{max}}(2) \\ \vdots & \vdots & \vdots & \ddots & \vdots \\ r^i_{T_{max},1}(t) & r^i_{T_{max},2}(t) & r^i_{T_{max},t+l}(t) & \cdots & r^i_{T_{max},T_{max}}(t) \end{pmatrix}. \quad (4)$$

By virtue of the introduced notation (2), (3), the diagonal elements of the matrix $R^i(T_{max})$ are equal to unity, i.e. $r^i_{T_{max},T_{max}}(t)$ for all $T = \{t: t = 1, \ldots, T_{max}\}$ and all $t$, and the remaining elements are in the range from -1 to +1 ($-1 \leq r_{ij} \leq 1$). This matrix (5) allows for the determination of moment $t$ in which the change of the rules of grouping and fixing the values of $x(t)$ occurred.

$$G_j^{total} = \sum_{t=1}^{T_{max}} r^i_{1,t+l}. \quad (5)$$

The integral indicator of the entire system of the $i$ - th enterprise.

$$G_i = \sum_{j=1}^{n} G_j^{total}. \quad (6)$$

IV. DESCRIPTION OF THE COMPANY SEZ AGAINST WHICH SANCTIONS ARE IMPOSED

The formation of residents of a special economic zone is a mechanism for the development of the economy of a constituent entity of the Russian Federation in its territory, aimed at improving the competitive characteristics of the Russian economy. This mechanism allows us to accelerate the pace of economic development of enterprises of a constituent entity of the Russian Federation by creating preferential taxation regimes.

In order to receive preferential tax treatment, an enterprise must go through the procedure of obtaining the status of a resident of the SEZ. The inclusion and maintenance of the list of residents of the SEZ is carried out by the federal executive body authorized by the Government of the Russian Federation. All activities of the SEZ are regulated by the federal law "On Special Economic Zones in the Russian Federation" of July 22, 2005, No. 116-FZ.

The enterprise, a resident of a special economic zone and used in modeling, is engaged in added-value wood processing with non-waste production and the sale of finished products abroad. The procurement of 1,300 thousand cubic meters plank timber in the North-Yenisei district is carried out. Delivery of the harvested raw materials is carried out on barges along the Yenisei River during the shipping period from June to September. Products from wood processing line are produced from harvested round timber: floorboards, glued beams, eurolining, etc. The production model takes into account acquired fixed assets: harvesting complexes, logging trucks, wood processing line, pellet equipment, boilers and main production equipment maintenance technology. The production process begins after the stage of procurement of a sufficient amount of raw materials in warehouses. Output is uniform throughout the year. From the first half of the first year of the period under consideration, the company has been undergoing the implementation of project management according to the PMBoK standard. By the beginning of the third year, equipment has been purchased to indouble the procurement and processing of raw materials. Activities, investment and organizational measures are carried out at the expense of the bank's credit facilities and the use of tax incentives for SEZ residents. The SEZ resident's modes of operation were calculated using the software for calculating the economic model of an enterprise engaged in harvesting and wood processing line of various wood species (larch, pine, aspen, birch, cedar) under certain market development scenarios and enterprise strategy [18].

Below is the structure of the simulated enterprise data.

- Project Owners. Pledge of the share of authorized capital in the bank. The structure of the shares of several owners.

- Parameters of products. Considered the size of semi-finished products, finished products, recyclable waste. Formed a balanced model of the movement of inventories in warehouses.

- Market research: the structure of the forest complex, domestic consumption, external consumption, analysis of the development of market trends.

- Plan of measures: stages of project implementation, including the holding of new investment and organizational measures.

- Marketing: product benefits accounted for; The analysis of target markets in Europe: Germany, France, South Korea, Japan; calculated market potential, supply structure.

- Plan of production: completed engineering process; The scheme of material production flows has been drawn up.

- Repair of fixed assets: the calculation of technical inspection and repair, overhaul.

- Cost analysis: the date of acquisition and entry of fixed assets, the characteristics of construction work, analysis of personnel decisions.

- Ecology and BC: the cost of job security and environmental safety of the facility.

- Project risks: political risks, risk of increased competition, risks of inadequate qualified personnel, risks of lack of demand for manufactured products, risks of inadequate qualifications of the company's management, risk of disruption to the start-up of production, risk of raw materials supply.

- Financial model: integral indicators of the project [15], asset analysis, liquidity analysis, analysis of business performance, analysis of the integral indicator, calculation of indicators of project sensitivity to risks, analysis of cost and project costs, alternative ways to implement the project, cash flow, profit and losses, balance, cost analysis of all stages of production, analysis of the turnover of work in progress, analysis of fixed assets, analysis of personnel costs, analysis of operating costs of machines and mechanisms.

Variants of modeling the activities of enterprise of an SEZ resident depending on the microeconomic parameters are presented in a separate paper [19].

## V. RESEARCH RESULTS

The activity of the special economic zone of the Krasnoyarsk Territory $\Omega$ and each enterprise of the resident $X^i$, is modeled on 9.6 million $x_N^i$ values monthly for 5 years in the author's software complex [18, 20-22]. The results of management by the regional authorities (subject of the Russian Federation), for all parameters $u(t)$ listed above, are shown in Fig. 1. The calculation of the effect of each measure $u(t)$ was calculated in turn. At the output, we have the number $y(t)$ value equal to the developed measures $u(t)$.

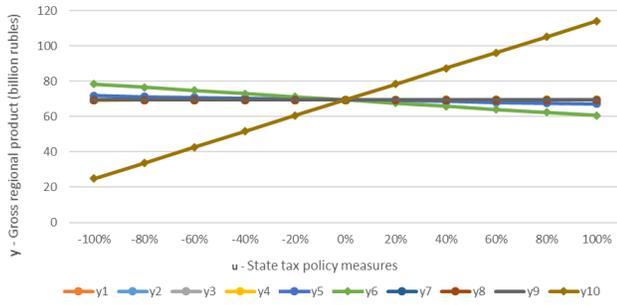

Fig. 1. Change of y from the investment policy regimes of the region of the Russian Federation

Greater efficiency for $y_6$ - Krasnoyarsk Territory GRP is achieved from $u_6$ - improving logistics conditions, for the example, paying subsidies for the reimbursement of a part of transportation costs for product delivery in the amount of 20 to 60%. Changes in $\Delta u_6$ - 20% of subsidies for the cost of transporting products will on average affect the annual $\Delta y_6$ - increase in the GRP of the Krasnoyarsk Territory in the amount of 0.2%. Fluctuations in the uncontrolled endogenous factor $v_{10}$ (dollar rate) by an average of 20% change the $y_{10}$ annual GRP by 2%.

According to the simulated situation, the inflow of resources to the enterprise is limited from 37th periods fig. 2.

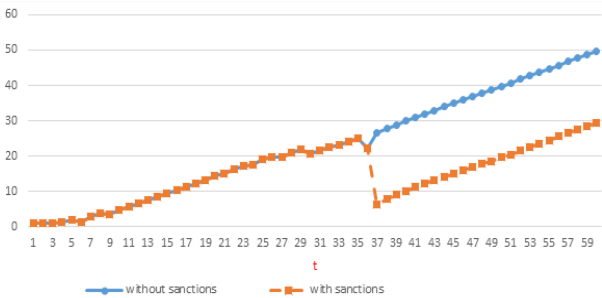

Fig. 2. The dynamics of the integral indicator in different modes of operation of the system

Fig. 2 shows that the calculation of the integral index $G_i(t)$ changes its values in the period when the sanctions began to affect the enterprise. In our case, this is the 37th period. The values of the integral index in the operation mode without sanctions are $G_i$ =1,369, in the mode of sanctioning $G_i$ =887.

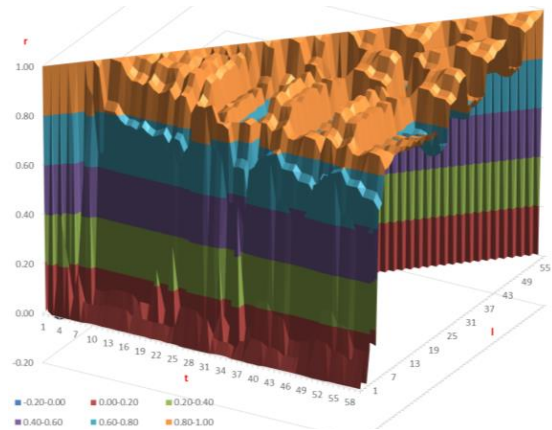

Fig. 3. Three-dimensional modeling of the integral indicator of the enterprise without sanctions

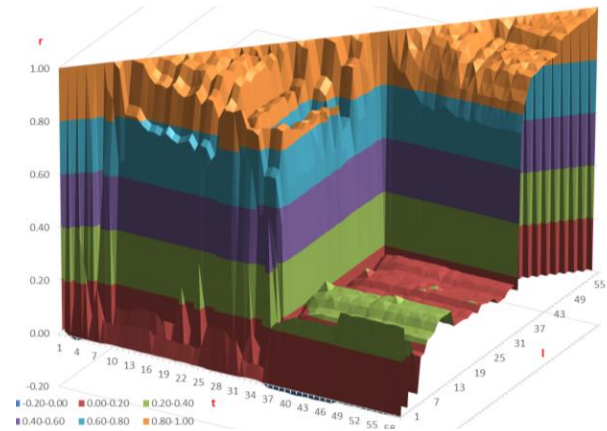

Fig. 4. Three-dimensional modeling of the integral indicator of the enterprise with sanctions

Calculation of $R^i(T_{max})$ by the formula 4 allows us to consider Fig. 2 in a three-dimensional plane. From fig. 3 it can be seen that the enterprise does not have changes in the accounting policy of income and expenses a dynamic system does not change the structure of accounting for its parameters characterizing $x(t)$. Fig. 4 shows the imposed sanctions in the 37th period and their further impact on the enterprise.

## VI. THE DISCUSSION OF THE RESULTS

The restriction of export sanctions by the company to 80% in the 37th period forces the company to look for other sales markets. The amount of revenue is not enough to finance its costs, so the company needs to raise borrowed funds with a higher percentage. For three months, it is possible to restore product sales to 20%, the rest of the products are shipped to the warehouse. In this situation, during the year the company is waiting for a technical default, as there are not enough financial resources. Within one year, the federal authorities will not have time to work out and implement supporting measures. In the formed conditions, the company is forced to reduce funding for production and administrative functions 5 times, selling off assets and closing projects.

The results of the experiment showed that for the Krasnoyarsk Territory the gross regional product (GRP) and gross domestic product (GDP) for 5 years from the work of the special economic zone with the imposed sanctions will decrease by 17.3% from 69 billion rubles to 57 billion rubles.

The obtained experimental data are identical to the current situation with the sanctions against RUSAL, introduced in April 2018 by the US Treasury Department. It imposed sanctions that reduced the company's exports by 80%, dropping shares by 46.9%. The remaining assets of the businessman fell by 18%, the state of the owner was reduced by half. The company announced news of a possible technical default. Within six months, RUSAL closed the city-forming enterprise in Karelia (Nadvoitsky Combine). Taking into account the resulting debt on loans and falling prices on the commodity market under the threat of closing the plant in Sweden (Kubal), in Khakassia (SAAZ) and other plants.

## VII. CONCLUSION

The activity of a special economic zone is described by a dynamic equation. Modeling showed an increase in the GRP of a constituent entity of the Russian Federation by 0.2-2% with the implementation of measures that reduce the logistical costs of SEZ residents by 20% and the market appreciation of the dollar by 20% each year for 5 years. The activity of each SEZ resident is modeled according to its individual parameters of assets, economic activity, goals and strategies. The sanctions regime, the export limit of 80%, to one of the residents of the special economic zone was applied. Formed integral indicators recorded simulated economic sanctions in the 37th period: $G_i$ =1,369 without sanctions and $G_i$ =887 taking into account sanctions. The socio-economic assessment of losses arising for the enterprise, the Krasnoyarsk Territory and the country was carried out. It has been proved that with powerful sanctions, which deprive an enterprise of 80% of its sources of funding, there is a loss of all economic advantages, market shares and other economic incentives to continue its work. A comparison with the actual situation allows us to conclude that the state is not in a position to implement compensatory measures that mitigate the impact of sanctions during the year.

The purpose of this work, indicated at the beginning of the article, has been achieved.